\documentclass[10pt]{amsart}

\usepackage[foot]{amsaddr}

\makeatletter
\renewcommand{\email}[2][]{%
  \ifx\emails\@empty\relax\else{\g@addto@macro\emails{,\space}}\fi%
  \@ifnotempty{#1}{\g@addto@macro\emails{\textrm{(#1)}\space}}%
  \g@addto@macro\emails{#2}%
}
\makeatother

\usepackage{lipsum}
\usepackage{amsmath,amsthm,amssymb,epsfig,hyperref,enumerate,array,tabularx,multirow}
\usepackage{graphicx,caption,natbib,fancyhdr}

\title[Dynamic Consistency]{Construction and Analysis of a Discrete Heat Equation Using Dynamic Consistency: The Meso-scale Limit}  

\date{}

\keywords{Heat equation; Dynamic consistency; Random walk; Asymptotic analysis; Continuum limit} 

\subjclass[2020]{35K05,; 39A14}

\theoremstyle{definition}
\newtheorem{definition}{Definition}[section]

\theoremstyle{comm}
\newtheorem{comm}{Comment}[section]

\begin{document}

\author{Ronald E.~Mickens$^a$}
\address[a]{Department of Physics, Clark Atlanta University, Atlanta, GA, 30314} 
\email[a]{rmickens@cau.edu}

\author{Talitha Washington$^b,^*$}
\address[b]{Department of Mathematical Sciences, Clark Atlanta University, Atlanta, GA, 30314} 
\address[$^{\ast}$]{Corresponding author}
\email[b]{twashington@cau.edu}

\begin{abstract}
We present and analyze a new derivation of the meso-level behavior of a discrete microscopic model of heat transfer. This construction is based on the principle of dynamic consistency.  Our work reproduces and corrects, when needed, all the major previous expressions which provide modifications to the standard heat PDE.  However, unlike earlier efforts,  we do not allow the microscopic level parameters to have zero limiting values.  We also give insight into the difficulties of constructing physically valid heat equations within the framework of the general mathematically inequivalent of difference and differential equations. 
\end{abstract}

\pagestyle{fancy}
\fancyhead[R]{ }

\maketitle

\section{Introduction}\label{Introduction}

The purpose of this paper is to examine the meso-scale limit of a discrete micro-level mathematical model constructed to represent simple heat transfer.  The creation of the microscopic model is based on the application of the principle of dynamic consistency (\cite{13,  14,  15}).  Under the appropriate mathematical assumptions, we are able to obtain the previous results of  \cite{11}, \cite{4}, and \cite{16} regarding the replacement of the standard heat equation
\begin{equation}
u_t=Du_{xx},\label{eq1}
\end{equation} by the generalization
\begin{equation}
\tau u_{tt}+u_t=Du_{xx},  \label{eq2}
\end{equation}
where $D$ is the temperature diffusion constant,  $\tau$ is a time-lag parameter,  and
\begin{equation*}
u_t = \frac{\partial u}{\partial t},  \quad u_{tt} = \frac{\partial^2u}{\partial t^2},  \quad u_x=\frac{\partial u}{\partial x},  \quad u_{xx}=\frac{\partial^2 u}{\partial x^2}. 
\end{equation*}
Note that $(x,t)$ are the one-dimensional space and time independent variables, and $u=u(x,t)$ is defined over appropriate intervals of $(x,t)$ with suitable boundary conditions and initial values. The need for a generalization of Eq. ~(\ref{eq1}) comes from the fact that the solutions of Eq.~(\ref{eq1}) transmit information at an infinite speed (\cite{2,5,6,7,8}),  a condition which violates the principle of causality (\cite{2,5,6,7,8}).

For convenience,  we now provide a resume of our major results:
\begin{itemize}
\item[(i)] Based on our formulation of the discrete mirco-level model, the drived first modified form of the standard heat PDE has the mathematical structure
\begin{equation} \label{eq3}
\epsilon_1 u_{tt} +u_t= Du_{xx}+(\epsilon_2D)u_{xxxx}.
\end{equation}
where $(\epsilon_1, \epsilon_2)$ are non-negative parameters.
\item[(ii)] Finite-order approximate models of discrete micro-level models can never be ``exact." By this, we mean that the meso-level equations will always require knowledge of certain properties of the solutions which in practice and/or theoretically can not be experimentally measured. 
\end{itemize}

This paper is organized as follows: Section 2 provides a brief historical summary of work on the standard heat PDE,  its major defect,  and previous attempts to resolve the associated issue. It also includes several of the mathematical models put forth to solve the causality problem and discusses why they all have failed in this quest. Section 3 gives several general comments on the modeling process and introduces the concept of dynamic consistency.  In Section 4, we use dynamic consistency to construct perhaps the simplest discrete mathematical model for heat conduction; it corresponds to a generalized random walk model.  We take this to be our micro-level model. In Section 5, we show that by using a Taylor series expansion, a meso-level hybrid difference-differential equation is obtained.  Special cases of this equation includes all previous cases of the modified standard heat equation.  Finally,  in Section 6,  we summarize our major conclusions and provide reasons why the usual methods for generalizing the heat PDE may never produce valid results.

A last comment.  We have not tried to reference every relevant publication on the subjects of this work.  Such references are generally well known by active researchers in these areas and are readily accessible.  This view is consistent with the fact that the current article is not a review paper.

\section{Preliminaries}\label{Preliminaries}

Elementary heat conduction phenomena have been modeled by for approximately 150 years by the linear partial differential given in Eq.~(\ref{eq1}).  However, this PDE has associated with it two difficulties.  First,  Eq.~(\ref{eq1}) allows for the transformation of information between two separate space points with infinite speed (\cite{8}).  Second, this PDE does not provide a valid or accurate model for many heat conduction problems arising in a broad range of scientific and engineering situations (\cite{6,7,8}).  As a result of the existence of these issues,  there has been a large research effort to resolve and/or explain these issues; see for example (\cite{2, 5,6,7,8}). 

In general,  the modifications to Eq.~(\ref{eq1}) have included adding a second-derivative term in time. The simplest such model is the Maxwell-Cattaneo representation (\cite{11, 4, 16}),
\begin{equation} \label{eq4}
\tau u_{tt} +u_t = Du_{xx} 
\end{equation}
where the constant,  positive parameter $\tau$ is given the interpretation of a delay or lag time. Note that this PDE has exactly the same mathematical structure as the damped wave equation (\cite{2,7,11}).  Note that since the dependent variable,  $u(x,t)$,  represents the temperature as measured with respect to the absolute Kelvin scale, then only solutions satisfying the condition
\begin{equation} \label{eq5}
u(x,t) \geq 0, \quad t\geq 0,
\end{equation}
are physically meaningful. Inspection of Eq.~(\ref{eq4}) shows that it is a hyperbolic PDE and,  as a consequence,  can produce negative valued solutions.  So this is a difficulty that must be overcome.  However,  Eq.~(\ref{eq4}) does resolve the issue of infinite speed for information transfer.  We expect this speed to be approximately
\begin{equation} \label{eq6}
c=\sqrt{\frac{D}{\tau}}.
\end{equation}
Other modifications of Eq.~(\ref{eq1}) produce the linear PDE
\begin{equation*}
\tau u_{tt} +u_t=Du_{xx}+D_1u_{xxt},
\end{equation*}
where $(\tau,  D,  D_1)$ are non-negative,  constant parameters; see,  for examples, the paper of \cite{7}.  Here,  the issues also include whether the data can actually be measured to obtain the appropriate initial and boundary conditions.

Our analysis of these modified heat equations leads us to the conclusion that these efforts have no fundamental groundings in an overall consistent physical theory.  Further, they are related to other efforts involving the examination of macroscopic limits of microscopic mathematical models (\cite{1,3}) and the validity of expansions in a (supposedly) small parameter (\cite{12}).  An accurate,  but colorful way of characterizing this situation is to look at it from the perspective of a ``wack-a-mole" situation where the things required to resolve one issue provides the opportunity for the creation of another unresolved issue.

\section{Modeling and dynamic consistency}\label{Modeling}

A mathematical model of a system is a symbolic representation incorporating some aspects of its properties and features.  These aspects concern properties and important features of immediate interest to the modeler and the relevant scientific community.  However,  it must not be forgotten that a mathematical model,  in general,  can never capture all aspects and details of a system.  This is because at a given moment, no realistic system has all of its properties known and as a consequence,  the unknown ones can not  be included in the model.  Further,  even if a particular feature or aspect of a system is known,  the currently available mathematical structures may not be able to incorporate them in an appropriate and viable model.

In the next Section we present arguments which form the bases for the construction of a discrete microscopic model of heat conduction in one-space dimension.  Our dependent variable,  $u$,  represents the temperature measured on the absolute Kelvin scale.  This model will be derived by using the principle of dynamic consistency (\cite{13,15}).  We now provide a concise definition of dynamic consistency and explain how it may be applied to the construction of mathematical models.

\begin{definition}
Consider two systems, $S_1$ and $S_2$.  Let $S_1$ have some property, $P$.  If $S_2$ also has property $P$, then $S_2$ is said to be dynamically consistent to $S_1$ with respect to property $P$.
\end{definition}

Observe that the systems can be essentially anything. An explicit example is where $S_1$ is a differential equation and $S_2$ is a discretization of $S_1$.  The case where $S_2$ is a finite-difference discretization has been extensively studied within the context of the ``nonstandard finite differences” methodology (\cite{13,15}).

Another example is where $S_1$ is a physical system and $S_2$ is a mathematical model of $S_1$. For this case, $P$ might include any or all of the following features (\cite{13}):
\begin{enumerate}[(i)]
\item the dependent variables are non-negative
\item the dependent variables are bounded
\item conservation laws exist
\item the dependent variables have special asymptotic behaviors.
\end{enumerate}
Note that the mathematical model may not have all of the properties $P$. This is a general consequence of the modeling process. Again, considering the situation where $S_1$ is a differential equation and $S_2$ its discretization, then numerical instabilities may occur in the equations of discretizations (\cite{13,15}).

Since the major goals of mathematical modeling are to provide an understanding of the system, help to make valid predictions of the evolution of the physical system, and otherwise provide insights into the fundamental nature of the system,  we must be thoughtful in our construction of such models.  In particular,  we must always be aware that the results of mathematical modeling depend strongly on the mathematical structures we use to create the models.  An insightful example is the modeling of single-population systems where discrete models give solution behaviors entirely different from those provided by differential equations (\cite{14}). 

Another difficulty in creating mathematical models is to take as exact particular functional forms which are based on experimental results valid only at low concentrations, or a limited range of temperatures, etc.  A prime illustration of this is Fick’s law for heat conduction (\cite{6,8}).

Finally,  we give a concise summary of the dynamic consistency methodology for constructing mathematical models of a system: 
\begin{enumerate}[(a)]
\item Understand and analyze the system as fully possible; 
\item from (a),  list the major properties and constraints that should be incorporated into the model;
\item using a particular mathematical structure (difference equations, differential equations, integral equations, etc.),  construct a mathematical model that includes as many as possible of the features and/or items listed in (b). 
\end{enumerate}
Note that different mathematical models will arise from the effort to include some, all, or a separate set of properties of the system and its associated constraints,  such as symmetries and conservation laws.  In general,  the application of the dynamic consistency methodology does not yield a unique mathematical model (\cite{9,13}],  and often this ambiguity can be used to our advantage.  It should be clear that the inclusion of an increasing number of dynamic, consistent properties should yield better mathematical models.  However,  the modeler may not,  in general, know all the important properties of a system.  This implies that it is unlikely that the modeling process will generate equations that are ``exact,” either mathematically or physically.

In the next section,  we apply the dynamic consistency methodology to the issue of constructing a discrete microscopic level model for simple heat conduction. 

\section{A discrete microscopic level model}\label{Discrete}

Before proceeding to construct the microscopic level model equation,  we need to provide some brief comments on several topics.

\begin{figure}
\centering
\includegraphics[scale=0.1]{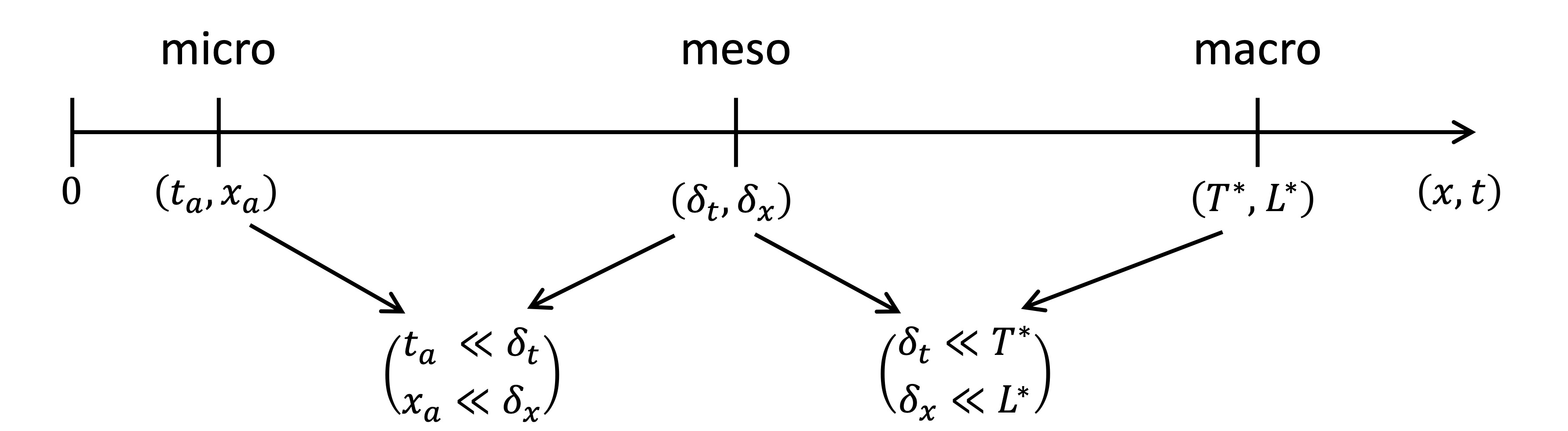}
\caption{Ordering of time and space scales}\label{fig1}
\end{figure}

First,  we assume that it makes sense to divide the physical universe into three separate relevant domains in terms of space-time characterization.  This is depicted in Fig. ~(\ref{fig1}) and Table~(\ref{tbl1}).  The smallest or lowest level is where the atomic structure of matter dominates our understanding of physical systems. The top level,  sometimes called the continuum limit, is where macro-scale phenomena is often observed,  measured,  and interpreted.  Generally,  at this level,  the graininess of the atomic aspects of matter does not play an important role. Finally, the meso-level corresponds to considering an a priori atomic system,  but averaged over its atomic distinctness to produce variables that have macroscopic meaning.

Second,  we associate with each of the three levels characteristic time and space physical scales; see Fig. ~(\ref{fig1}). and Table~(\ref{tbl1}).  Thus, $(t_a, x_a)$ are constant parameters associated with the evolution of the system dynamics at the micro-level; $(\delta_t,\delta_x)$ are the constant parameters associated with physical phenomena at the meso-level; and $(T^*, L^*)$ are the time and length scales for what occurs at the macro-level.  Relations between their magnitudes are indicated in Fig.~(\ref{fig1}) and Table~(\ref{tbl1}).

%

\begin{table}
\caption{Definition of time and space related variables}\label{tbl1}
\begin{tabular}{| l | l | l}
$u_s^r$  & $u_m^k$ & $u(x,t)$  \\ 
$t_r=t_a\cdot r$& $t_k = \delta_t\cdot k$ & $(x,t) \textrm{ continuous}$\\
$x_s = x_a\cdot s$ & $x_m = \delta_x \cdot m$ & \\
$r=0,1,2,\dots$ & $k=0,1,2, \dots$ & \\
$s= \textrm{ integers}$ & $m=\textrm{ integers}$ & \\
\end{tabular}
\end{table}

Our main purpose is to construct a discrete space-time model for the micro-level and then extend it to the meso-level.  Next,  we examine the meso-level model and interpret it as a modified or corrected version of what is usually considered to be a mathematical model for the continuum limit.

The space and time variables at the micro-level are discrete and are defined in Fig.~(\ref{fig1}) and Table~(\ref{tbl1}).  We show that the derived meso-level evolution equation has to be considered a hybrid differential-difference equation. This we explain below.  Note that $t_r$ and $x_s$ are discrete time and space variables.  The discrete meso-scale variables and their relationship to the continuum variables $(t,x)$ will also be clarified later.

To proceed in the derivative of a discrete micro-level model for heat conduction in one-dimension, we assume that our mathematical model satisfies the following conditions:
\begin{enumerate}[(A)]
\item be a linear difference equation,
\item satisfy the constraint that information can only travel with a finite speed in the discrete space-time;
\item be invariant under the parity transformation; i.e., the discrete equation should not change its form under the interchange
\begin{equation}\label{eq7}
(s+i) \ \longleftrightarrow \ (s-i),\\
s = \textrm{ integer}, \ i=1, 2, 3, \dots.
\end{equation}
\item not be invariant under time inversion;
\item satisfy a non-negativity condition for its solutions;
\item possess the property of measurability.\\
\end{enumerate}

The following is a brief discussion of the implications of these requirements:
\begin{enumerate}
\item[(A$_1$)] Assuming that the micro-level equation is linear provides a great simplification in that there are many more possible nonlinear equations and, generally, we have no clue as to which one(s) to select.
\item[(B$_1$)]  The micro-level is characterized by the time scale to $t_a$ and the space scale to $x_a$. Therefore, the speed of information flow,  $c$,  will be
\begin{equation} \label{eq8}
c = \frac{x_a}{t_a}
\end{equation}
\item[(C$_1$)]  For a discrete equation the condition given by Eq.~(\ref{eq7}) ensures parity conservation (\cite{13}).
\item[(D$_1$)]  The simplest way to have this condition hold is to have the difference equation depend only on the discrete times $t_r$ and $t_{r+1}$. Physically, this requirement implies physical dissipation.
\item[(E$_1$)]  Since our dependent variable corresponds to temperature on the absolution Kelvin scale, it must be always non-negative (\cite{6,11}).
\item[(F$_1$)]  By having the property of measurability, we mean that the model should not contain any parameter and/or dependent variables that can not be readily measured using standard experimental techniques.
\end{enumerate}

Let $u_s^r$ be the temperature at discrete time $t_r$ and location $x_s$.  The simplest difference equation that satisfies all of the above requirements is given by the expression
\begin{equation}
\begin{cases}\label{eq9}
u_s^{r+1} = &pu_{s+1}^r+(1-2p)u_s^r+pu_{s-1}^r\\
&0<p \leq \frac{1}{2}.
\end{cases}
\end{equation}
In the work to follow, we specialize to the case $p=1/3$ and Eq.~(\ref{eq9}) reduces to 
\begin{equation} \label{eq10}
u_s^{r+1} = \left(\frac{1}{3}\right)\left(u_{s+1}^r+u_s^r+u_{s-1}^r\right).
\end{equation}
This special value of $p$ does not change any of our conclusion.

It should be indicated that Eq.~(\ref{eq10}) is a partial difference equation, first-order in the discrete-time and second-order in the discrete space variables.  Furthermore,  Eq.~(\ref{eq10}) can be solved exactly by a number of standard techniques (\cite{14}).  We will not provide such solutions since their explicit expressions play no role in the purpose of this article,  namely,  the derivation of the meso-level equation and how it should be interpreted.  However,  we do assume that there exists an analytic function, $U(x,t)$, such that
\begin{equation} \label{eq11}
\begin{cases}
u_s^r = U(x_s, t^r),\\
s=\textrm{ integer }; \quad r=(0,2,2,\dots); \quad x_a \textrm{ and } t_a \textrm{ are fixed.}
\end{cases}
\end{equation}
Note that Eq.~(\ref{eq10}) can be rewritten to the form
\begin{equation} \label{eq12}
u_s^{r+1}-u_s^r = \left(\frac{1}{3}\right)\left(u_{s+1}^r-2u_s^r+u_{s-1}^r\right),
\end{equation}
and from this we can obtain the expression
\begin{equation} \label{eq13}
\frac{u_s^{r+1}-u_s^r}{t_a} = D\left(\frac{u_{s+1}^r-2u_s^r+u_{s-1}^r}{x_a^2}\right),
\end{equation}
where
\begin{equation} \label{eq14}
D=\left(\frac{1}{3}\right)\left(\frac{x_a^2}{t_a}\right).
\end{equation}
If $D$ is assumed constant and if the limits
\begin{equation}\label{eq15}
t_a \rightarrow 0,  \quad x_a \rightarrow 0,  \quad D = \textrm{ constant},
\end{equation}
then Eq.~(\ref{eq13}) becomes the following partial differential equation
\begin{equation} \label{eq16}
u_t=Du_{xx}.
\end{equation}

\begin{comm}
In more detail, the following are the correct limits:
\begin{equation}
\begin{cases}
t_a \rightarrow 0, x_a \rightarrow 0; \quad  r \to \infty,  \quad s \to \infty;\\
\left(\frac{1}{3}\right)\left(\frac{x_a^2}{t_a}\right) = D = \textrm{ constant};\\
t_ar = t = \textrm{ fixed},  \quad x_as = x = \textrm{ fixed.}
\end{cases}
\end{equation}
\end{comm}

Observe that Eq.~(\ref{eq16}) is the standard heat equation (\cite{8}).  Also, the procedure used to derive this result opens up a number of important issues related to the macroscopic limits of microscopic models (\cite{1,3,9}).  As we shall see in the next section,  an alternative analysis and interpretation is possible.

\section{Macro-level approximation}

Consider Eq.~(\ref{eq10}) and as stated above,  assume that there exists a function, $U(x,t)$,  with appropriate analytical properties in $x$ and $t$,  such that
\begin{equation} \label{eq18}
u_s^r=U(x_x,t_r).
\end{equation}
This implies that $U(x,t)$, for fixed $\overline{x}$ and $\overline{t}$,  satisfies the same recurrence/different equation as $u_s^r$, i.e., 
\begin{equation} \label{eq19}
U(x,t+\overline{t})=\left(\frac{1}{3}\right)\left[U(x+\overline{x},t)+U(x,t)+U(x-\overline{x},t)\right]
\end{equation}
Thus, at the meso-level where $(\delta_x,\delta_t)$ are the length and time (physical) scales,  we may identify $(\overline{x},\overline{t})$, respectively, with $(\delta_x,\delta_t)$,  and obtain
\begin{equation}\label{eq20}
U(x,t+\delta_t)=\left(\frac{1}{3}\right)\left[U(x+\delta_x,t)+U(x,t)+U(x-\delta_x,t)\right]
\end{equation}
To proceed,  carry out Taylor series expansions of the individual terms in Eq.~(\ref{eq20}) and retain only terms up to $O(\delta_t^3)+O(\delta_x^6)$.  Doing this gives, after division by $\delta_t$, the expression
\begin{equation}\label{eq21}
U_t+\left(\frac{\delta_t}{2}\right)U_{tt}+O(\delta_t^2) = DU_{xx}+\left[\left(\frac{D}{36}\right)\delta_x^2\right]U_{xxxx}+O(D \delta_x^4),
\end{equation}
where
\begin{equation}\label{eq22}
D=\frac{\delta_x^2}{3\delta_t}.
\end{equation}
Note by assuming that $D$ is fixed in value,  we have
\begin{equation}\label{eq23}
\delta_t=O(\delta_x^2).
\end{equation}
Consequently,  to terms of $O(1)$ and $O(\delta_t)+O(\delta_x^2)$, the evolution equations for $U(x,t)$ are
\begin{align}
O(1) &: \overline{U}_t=D\overline{U}_{xx},\label{eq24} \\
O(1)+\mathcal{O}(\delta_t)+O(\delta_x^2) &: \left(\frac{\delta_t}{2}\right)\overline{U}_{tt}+\overline{U}_t=D\overline{U}_{xx}+\left[\left(\frac{D}{36}\right)\delta_x^2\right]\overline{U}_{xxxx}. \label{eq25}
\end{align}
We have placed bars over the $U$’s to indicate that these equations only yield approximations to the exact PDE for $U(x,t)$.

Note that the lowest order expansion,  i.e.,  retaining terms to $O(1)$, the evolution equation is the standard heat conduction PDE.  Since
\begin{equation}\label{eq26}
\partial_t=D\partial_{xx}+O(\delta_t)+O(\delta_x^2),
\end{equation}
we have
\begin{align}\label{eq27}
\partial_{xxxx} = &\partial_{xx}\cdot \partial_{xx}\cr
= &\partial_{xx}\left[\left(\frac{1}{D}\right)\partial_t\right]+O(\delta_t)+O(\delta_x^2).
\end{align}
Substituting this latter result into the second term on the right side of Eq.~(\ref{eq25}) gives
\begin{equation} \label{eq28}
\left(\frac{\delta_t}{2}\right)\overline{U}_{tt}+\overline{U}_t=D\overline{U}_xx+\left(\frac{\delta_x^2D}{36}\right)\overline{U}_{xxt}.
\end{equation}
This PDE is the same PDE that was derived by \cite{7} under different assumptions.

\section{Discussion}

In this work, we began by assuming that we have the exact micro-level evolution equation for simple heat conduction,  given by Eq.~(\ref{eq10}).  Assuming that the solution of this discrete equation can be represented as the values of an analytic function,  $U(x,t)$,  evaluated on the lattice of discrete space and time points,  we used Taylor series expansion on this expression at the meso-level to obtain what we will call modified hybrid-PDEs,  see for example Eq.~(\ref{eq28}).  Dropping the bar over $\overline{U}$,  the following two lowest order PDEs were found
\begin{equation} \label{eq29} 
U_t=DU_{xx},   
\end{equation}
\begin{subequations}
\vspace{-.3in}
\begin{align}
\tau U_{tt}+U_t = &DU_{xx}+D_1U_{xxt},  \label{eq30a}\\
\tau U_{tt}+U_t = &DU_{xx}+D_2U_{xxxx},  \label{eq30b}
\end{align}
\end{subequations}
where in this representation, we have expressed them, especially the Eqs.~(\ref{eq30a}, \ref{eq30b}),  in a general way suitable for use as phenomenological equations, where $(\tau,D,D_1)$ are to be determined from experiment. In principle, as indicated in the previous section, these parameters are determined by the time and length scales at the meso-level.

Observe that Eq.~(\ref{eq29}) is the standard heat equation. Within our interpretation, this PDE should be considered the macroscopic limiting equation of the microscopic Eq.~(\ref{eq10}).  We do not designate it by the name continuum limit (\cite{1,3}), since no micro-level or meso-level parameters are required to go to zero.  Also, note that within the context of this work, relationships should exist between the space and time scales at the micro- and meso-levels, e.g.,
\begin{equation} \label{eq31}
\delta_x=N_1x_a, \quad \delta_t=N_2t_a,
\end{equation}
where $N_1$ and $N_2$ are large integers.

Since neither the micro- or meso-levels parameters go to zero, we expect the small parameters $(\epsilon_1, \epsilon_2)$, where
\begin{equation}\label{eq32}
\epsilon_2=\frac{\delta_t}{T^*}, \quad \epsilon = \frac{\delta_x}{L^*},
\end{equation}
to play important roles in the expression of the first modification to the macroscopic equation.  To show this, carry out the scaling of independent variables in Eq.~(\ref{eq28}) as follows:
\begin{equation}\label{eq33}
t= T^*\overline{t}, \quad x=L^*\overline{x},
\end{equation}
where $(L^*, T^*)$ are the macro-level scales,  and $(\overline{x}, \overline{t})$ are the dimensionless independent variables.  Since the PDE is linear,  we need to make no explicit scaling for $U$. Substituting Eq.~(\ref{eq31}) into Eq.~(\ref{eq28}) and multiplying each term by $T^*$ gives
\begin{equation} \label{eq34}
\left[ \left( \frac{1}{2} \right) \left( \frac{\delta_t}{T^*}\right) \right] \frac{\partial^2 U}{\partial \bar{t}^2}+\frac{\partial U}{\partial \bar{t}}
=\left( \frac{DT^*}{(L^*)^2}\right)\frac{\partial^2U}{\partial \bar{x}^2}+\left( \frac{\delta_x^2}{36}\right)\left[ \frac{T^*D}{(L^*)^4} \right]\ \frac{\partial^4U}{\partial \bar{x}^4}.
\end{equation}
Keeping in mind the definition of $D$ given in Eq.~(\ref{eq22}), we define $\overline{D}$,  $\epsilon_1$, and $\epsilon_2$ as follows
\begin{equation} \label{eq35}
\overline{D}=\frac{DT^*}{(L^*)^2}, \quad \epsilon_1=\frac{1}{2}\left( \frac{\delta_t}{T^*} \right),  \quad \epsilon_2=\left( \frac{\delta_x}{L^*}\right)^2,
\end{equation}
where it should be noted that these three parameters are dimensionless,  i.e.,  pure numbers.  Also, it is expected that $\epsilon_1$ and $\epsilon_2$ are small,  i.e.,
\begin{equation}\label{eq36}
0<\epsilon_1<<1, \quad 0 < \epsilon_2<<1.
\end{equation}
With this notation,  Eq.~(\ref{eq28}) in dimensionless form is
\begin{equation} \label{eq37}
\epsilon_1 U_{\bar{t} \,\bar{t}}+U_{\bar{t}}=\overline{D}U_{\bar{x}\bar{x}}+\epsilon_2\overline{D}U_{\bar{x}\bar{x}\bar{x}\bar{x}}.
\end{equation}
From the fact that $\delta_t=O(\delta_x^2)$, we can expect an infinite hierarchy of modified PDEs,  starting from
\begin{equation} \label{eq38}
\begin{cases}
O(1): U_t=DU_{xx},\\
O(\delta_t)+O(\delta_x^2):  \textrm{Eq.~(\ref{eq37})}
\end{cases}
\end{equation}


\begin{figure}
\centering
\includegraphics[scale=0.1]{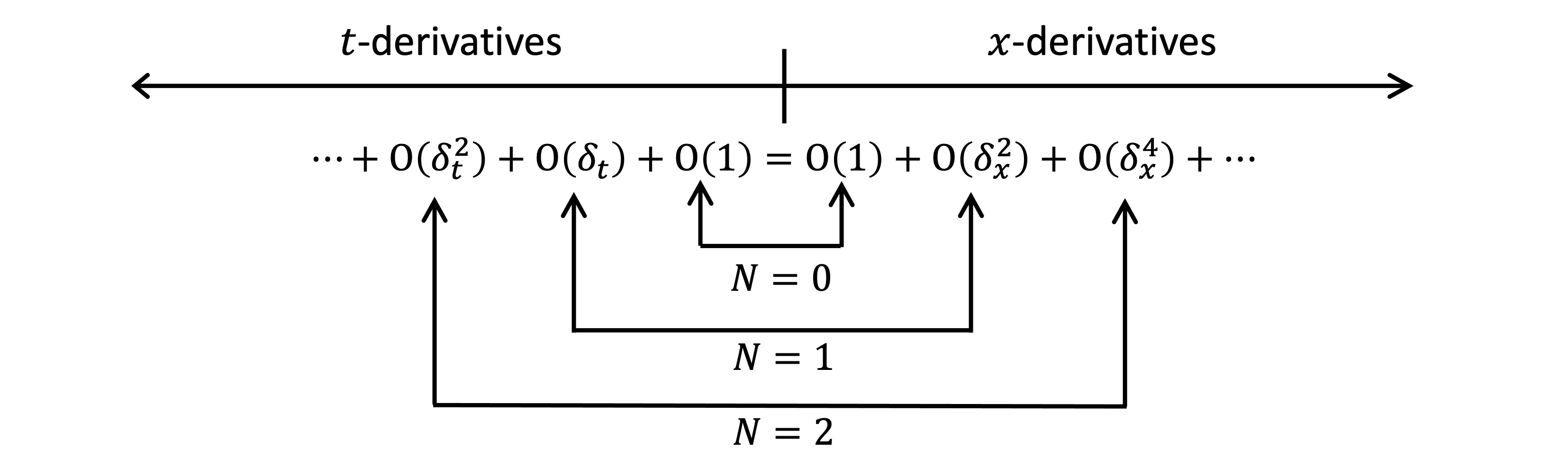}
\caption{Hiearchy of modified PDE.}\label{fig2}
\end{figure}

Keep in mind the fact the $N$-th modified PDE contains terms from all the previous equations. The $N=0$ and $N=1$ cases are the two expressions presented in Eq.~(\ref{eq38}).

Observe that at the $N$-th level of the modified PDEs hierarchy,  the time derivatives up to order $(N+1)$ appear.  This means also that spatial derivatives to even order $2(N+1)$ must occur.  Thus,  for $N=1$,  the time derivatives $U_t$ and $U_{tt}$ appear only along with the spatial derivatives, $U_{xx}$ and $U_{xxxx}$ or $U_{xxt})$.  However,  for most if not all one-dimension space systems,  the initial condition,  $U(x,0)$,  is known,  along with related boundary conditions, but $U_t(x,0)$ or higher time derivative are not known, it follows that none of the modified PDEs is measurable,  i.e.,  their mathematical solutions require information that can not be physically determined.  An important implication is that this way or methodology of investigating heat condition is invalid or at least not productive.  One way to view this situation is to remember that if we begin with a discrete space and time model, then the mathematical structure currently used corresponds to difference equations.  However,  it is well known that difference equations do not model exactly differential equations; see the books by (\cite{14,15}).  In particular,  the following relation holds
\begin{equation} \label{eq39}
\textrm{ difference equation } = \textrm{ infinite-order differential equation.}
\end{equation}
This can easily be seen by considering the following argument:
\begin{itemize}
\item[(1)] Define the (difference) shift and derivative operators by the following relations
\begin{equation}\label{eq40}
E_af(x) \equiv f(x+a), \quad Df(x)\equiv \frac{df(x)}{dx}.
\end{equation}
\item[(2)] The Taylor expansion of $f(x+a)$ is 
\begin{align}\label{eq41}
f(x+a)&=\overset{\infty}{\underset{k=1}{\sum}} \left(\frac{a^k}{k!}\right)D^kf(x)\cr
&=e^{aD}f(x).
\end{align}
\item[(3)] From the definition of the shift operator,  it follows that the operators satisfy
\begin{equation} \label{eq42}
E_a=e^{aD}
\end{equation}
\item[(4)] The difference operator is defined by the relation
\begin{equation} \label{eq43}
\Delta_a=E_a-1.
\end{equation}
Comparison of Eq.~(\ref{eq42}) and Eq.~(\ref{eq43}) gives
\begin{equation} \label{eq44}
\Delta_a=e^{aD}-1.
\end{equation}
\item[(5)] Since $\Delta_a/a$ and $D$,  are the discrete and continuous forms of the derivative, Eq.~(\ref{eq44}) shows that they nonlinearly related and any finite order difference equation corresponds to an infinite order differential equation.  Also,  if we solve Eq.~(\ref{eq44}) for $D$, we obtain
\begin{align}
 D &\equiv \left(\frac{1}{a}\right) \ln (1+\Delta_a)\\
&= \left(\frac{1}{a}\right)\left[\Delta_a-\frac{\Delta_a^2}{2}+\frac{\Delta_a^3}{3}-\frac{\Delta_a^4}{4}+\cdots\right].
\end{align}
Consequently,  it can be concluded that in general a finite order differential equation is equivalent to an infinite order difference equation.
\end{itemize}

In summary,  we have used the dynamic consistency methodology (\cite{13,14}) to construct a discrete space and time microscopic mathematical model for one-space dimension heat conduction.  To obtain solutions to this difference equation, only initial and/or boundary values are required to be known. After some manipulations,  we were able to derive an infinite number of mesoscopic PDEs.  These PDEs form a hierarchy of equations,  each of which modifies the equations on previous levels; see Fig. ~(\ref{fig2})  However, the modified equations at the first and higher levels generally require knowledge of time derivatives for which experimental data is either difficult or impossible to measure.  Our major conclusion is that the procedures used to construct and justify previous generalizations of the basic heat conduction equation,  $u_t=Du_{xx}$,  are not valid and that other mathematical structures must be discovered or created to resolve issues such as infinite speed of propagation for the transfer of information (\cite{2,4,5,6,7,8,11,16}).  We are currently investigating this problem.




%





\section{Bibliography}\label{Bibliography}


\bibliographystyle{cas-model2-names}

\bibliography{}



\end{document}